\documentclass[a4paper,12pt,reqno]{amsart}
\usepackage{babel}
\usepackage[utf8]{inputenc}
\usepackage[T1]{fontenc}
\usepackage{amsmath,amsthm,amssymb}
\usepackage{fullpage}
\usepackage{mathtools}
\usepackage{microtype}
\AtBeginDocument{\sbox0{$x$}} 
\usepackage[parfill]{parskip}
\usepackage{embedfile}
\usepackage[shortlabels]{enumitem}
\setlist[enumerate,1]{label={\upshape(\roman*)}}
\usepackage[
  pdftitle={Groups satisfying Gaschütz's complement theorem},
  pdfauthor={Stefanos Aivazidis, Benjamin Sambale},
  pdfstartview={FitH}
]{hyperref}
\AtBeginDocument{%
}

\newtheorem{prevtheorem}{Theorem}

\newtheorem*{namedtheorem}{Theorem}

\newtheorem{theorem}{Theorem}[section]
\newtheorem{corollary}[theorem]{Corollary}

\newcommand{\card}[1]{\left\lvert #1\right\rvert}
\newcommand{\gensub}[1]{\left\langle #1\right\rangle}
\newcommand{\gasch}[1]{\mathcal G(#1)}

\DeclareMathOperator{\zent}{\mathrm{Z}}
\DeclareMathOperator{\aut}{\mathrm{Aut}}
\DeclareMathOperator{\inn}{\mathrm{Inn}}
\DeclareMathOperator{\hol}{\mathrm{Hol}}

\makeatletter
\@removefromreset{equation}{section}
\makeatother

\newenvironment{proofof}[1][Proof]
  {\begin{proof}[{\upshape\bfseries #1}]}
  {\end{proof}}

\newenvironment{proofofref}[1]
  {\begin{proofof}[Proof of \autoref{#1}]}
  {\end{proofof}}

\title{Groups satisfying Gaschütz's complement theorem}
\author{Stefanos Aivazidis}
\address{Athens, Greece}
\email{\href{mailto:stefanosaivazidis@gmail.com}{stefanosaivazidis@gmail.com}}
\author{Benjamin Sambale}
\address{Institut für Algebra, Zahlentheorie und Diskrete Mathematik\\
Leibniz Universit\"at Hannover\\
Welfengarten 1\\
30167 Hannover\\
Germany}
\email{\href{mailto:sambale@math.uni-hannover.de}{sambale@math.uni-hannover.de}}
\date{\today}

\begin{document}
\embedfile[desc={GAP code deciding Gaschütz's theorem},
           mimetype={text/plain}]{gaschuetz.g}
\begin{abstract}
A classical theorem of Gaschütz states that an abelian normal
subgroup $N$ of a finite group $G$ has a complement in $G$ whenever
it has a complement in some subgroup $H$ with $N \leq H \leq G$ and
$\gcd(\card{N}, [G:H]) = 1$. We characterize the finite groups $N$
for which this theorem remains valid without the abelian hypothesis:
this is the case if and only if $\zent(N) \cap N' = 1$ and $N$ has a
complement in a quotient of the holomorph $\hol(N)$ by a diagonal copy of $N'$. 
This completes a previous paper by the second author.
\end{abstract}

\maketitle

\section{Introduction}\label{sec:introduction}

All groups considered in this paper are finite. A subgroup $C$ of a
group $G$ is a \emph{complement} of a normal subgroup
$N \unlhd G$ if $G = NC$ and $N \cap C = 1$. Deciding
whether a normal subgroup has a complement (that is, whether $G$
splits over $N$) is a classical problem of group theory.
For abelian $N$, a well-known theorem of Gaschütz reduces this
question to a subgroup of coprime index.

\begin{namedtheorem}[Gaschütz \cite{Gaschuetz1952}]
Let $N$ be an abelian normal subgroup of a group $G$, and let
$N \leq H \leq G$ with $\gcd(\card{N}, [G:H]) = 1$. Then $N$ has a
complement in $G$ if and only if $N$ has a complement in $H$.
\end{namedtheorem}

This result was extended to groups $N$ with abelian Sylow subgroups by the second author in \cite[Theorem~2]{Sambale2023}.\footnote{Note that the \href{https://arxiv.org/abs/2303.00254}{arXiv version} of this paper has a different numbering for the theorems.}
Following this paper, we say that
\emph{Gaschütz's theorem holds for} $N$, and write
$\gasch{N}$, if the theorem holds without conditions on $N$: whenever $N \unlhd G$,
$N \leq H \leq G$, $\gcd(\card{N}, [G:H]) = 1$, and $N$ has a
complement in $H$, then $N$ has a complement in $G$. 
By \cite[Theorem~12]{Sambale2023}, $\gasch{N}$ implies
$\zent(N) \cap N' = 1$, where $N'$ denotes the commutator subgroup of $N$. 
In particular, $\gasch{N}$ fails for all non-abelian nilpotent groups. 
For metabelian groups $N$, $\gasch{N}$ is in fact equivalent to $\zent(N) \cap N' = 1$ (\cite[Corollary~14]{Sambale2023}).

As indicated in \cite[Introduction]{Sambale2023}, Gaschütz and Brandis both expressed their wish to find a complete intrinsic characterization of $\gasch{N}$. This is the goal of the present paper.
For $x\in N$, let $\gamma_x$ be the inner automorphism of $N$ induced by $x$, i.e., $\gamma_x(y)\coloneqq xyx^{-1}$ for $y\in N$. Let $\hol(N) \coloneqq N\rtimes\aut(N)$ be the holomorph of $N$. 

\begin{prevtheorem}\label{thm:main}
Let $N$ be a group and
\[
D \coloneqq \bigl\{(d^{-1}, \gamma_d) : d \in N'\bigr\} \unlhd \hol(N).
\]
Then Gaschütz's theorem holds for $N$ if and only if
$\zent(N) \cap N' = 1$ and $N$ has a complement in
$\hol(N)/D$, where $N$ is identified with its
natural image in this quotient.
\end{prevtheorem}

For abelian $N$, the two
conditions hold trivially, in accordance with
Gaschütz's theorem. 
For perfect groups $N$, we recover \cite[Theorem~15]{Sambale2023}: $\gasch{N}$ holds if and only if $\zent(N)=1$ and $\inn(N)$ has a complement in $\aut(N)$.  

The proof consists of two reductions of independent interest. In
\autoref{sec:derived-subgroup} we show that $\gasch{N}$ is
equivalent to the condition that every extension $B$ of $N$ which
splits modulo the commutator subgroup $N'$ already splits over $N$
(\autoref{thm:derived}). The proof turns any extension violating
this condition into a Gaschütz counterexample
$N \leq H \unlhd G$ with $[G:H] = p$, for any prime
$p \nmid \card{B}$. 
In \autoref{sec:finite-extension} we prove
that, under the necessary condition $\zent(N) \cap N' = 1$, every
extension of $N$ which splits modulo $N'$ is a fiber product of the single
extension $\hol(N)/D$ of
$\aut(N)/\gamma(N')$ (\autoref{thm:universal}). 

In the last section, we give an explicit example by answering \cite[Problem~18]{Sambale2023}:
$\gasch{N}$ does not hold for $N\cong (C_3^2 \rtimes Q_8) \times C_2$.

\section{Splitting modulo the commutator subgroup}\label{sec:derived-subgroup}

The property $\gasch{N}$ quantifies over pairs of groups
$N \leq H \leq G$. Our first reduction replaces these pairs by single
extensions of $N$ and shows that only the splitting behavior modulo
the commutator subgroup matters. This naturally extends Gaschütz's original theorem.

\begin{theorem}\label{thm:derived}
For every group $N$, the following conditions are equivalent.
\begin{enumerate}
\item\label{item:gaschuetz-property} $\gasch{N}$ holds.
\item\label{item:derived-splitting} For every group $B$ with
      $N \unlhd B$, if $N/N'$ has a complement in $B/N'$,
      then $N$ has a complement in $B$.
\end{enumerate}
Moreover, if \autoref{item:derived-splitting} fails for $B$, then every
prime $p \nmid \card{B}$ gives groups
$N \leq H \unlhd G$ such that $G/H \cong C_p$ and $N$ has
a complement in $H$ but not in $G$.
\end{theorem}

\begin{proofof}
Suppose first that \autoref{item:derived-splitting} holds, and let
$N \leq H \leq G$ satisfy the hypotheses defining $\gasch{N}$, say
with a complement $K$ of $N$ in $H$. Since $N'$ is characteristic in
$N \unlhd G$, we have $N' \unlhd G$, and $N/N'$ is an
abelian normal subgroup of $G/N'$. Dedekind's modular law gives
$KN' \cap N = (K \cap N)N' = N'$, and $(KN')N = KN = H$; hence
$KN'/N'$ is a complement of $N/N'$ in $H/N'$. Since $\card{N/N'}$
divides $\card{N}$ and $[G/N' : H/N'] = [G:H]$, these two quantities
are coprime, so Gaschütz's theorem gives a complement of $N/N'$ in
$G/N'$. Applying \autoref{item:derived-splitting} with $B = G$ shows that
$N$ has a complement in $G$.

Suppose conversely that \autoref{item:derived-splitting} fails. 
We construct a witness $G$ for the failure of $\gasch{N}$, similar to \cite[proof of Theorem~12]{Sambale2023}.
Choose
$N \unlhd B$ such that $N/N'$ has a complement in $B/N'$ but
$N$ has no complement in $B$, and put $R \coloneqq B/N$. Let
$\pi \colon B \to R$ be the quotient map, choose a prime
$p \nmid \card{B}$, and let
\[
P \coloneqq \{(b_1, \ldots, b_p) \in B^p :
       \pi(b_1) = \cdots = \pi(b_p)\}
\]
be the $p$-fold fiber product of $\pi$ with itself. The cyclic
shift $\sigma$ of $P$, given by
$(b_1, \ldots, b_p) \mapsto (b_p, b_1, \ldots, b_{p-1})$, is an
automorphism of order $p$; put $Q \coloneqq P \rtimes \gensub{\sigma}$. The homomorphism
$P \to R$ sending each tuple to its common value
$\pi(b_1) = \cdots = \pi(b_p)$ is invariant under $\sigma$, so it
extends to a homomorphism $\tau \colon Q \to R$ with $\tau(\sigma) = 1$.
Form the fiber products
\[
G \coloneqq B \times_R Q
  = \{(b, q) \in B \times Q : \pi(b) = \tau(q)\},
\qquad
H \coloneqq B \times_R P = G \cap (B \times P),
\]
and let $\rho \colon G \to Q$, $(b, q) \mapsto q$, be the projection
onto the second coordinate. Given $q \in Q$, the surjectivity of
$\pi$ provides some $b \in B$ with $\pi(b) = \tau(q)$; hence $\rho$ is
surjective. We identify $N$ with
$\ker(\rho)=N \times 1 \unlhd G$. Since $H = \rho^{-1}(P)$ and
$Q/P \cong C_p$, we obtain $H \unlhd G$ and
$G/H \cong C_p$. Note that $\gcd(\card{N}, p) = 1$ because
$\card{N}$ divides $\card{B}$.

\emph{Step 1: $N$ has a complement in $H$.}

Let $\lambda \colon P \to B$ be the projection onto the first coordinate
and put
\[
K \coloneqq \{(\lambda(x), x) : x \in P\} \leq B \times P.
\]
Every $x \in P$ satisfies $\pi(\lambda(x)) = \tau(x)$, so $K \leq H$. If
$(\lambda(x), x) \in N \times 1$, then $x = 1$ and hence $\lambda(x) = 1$; thus
$K \cap N = 1$. Moreover, if $(b, x) \in H$, then
$\pi(b) = \tau(x) = \pi(\lambda(x))$, so $n \coloneqq b\lambda(x)^{-1}$ lies in
$\ker\pi = N$ and $(b, x) = (n, 1)(\lambda(x), x) \in NK$. Hence $H = NK$,
and $K$ is a complement of $N$ in $H$.

Aiming at a contradiction, we assume from now on that $N$ has a
complement $M$ in $G$.

\emph{Step 2: there is a homomorphism $\theta \colon Q \to B$ with
$\pi\theta = \tau$.}

By assumption,
$M \cap \ker\rho = M \cap (N \times 1) = 1$ and
$\rho(M) = \rho\bigl(M(N \times 1)\bigr) = \rho(G) = Q$, so $\rho$
restricts to an isomorphism $M \to Q$. Composing its inverse
$Q \to M$ with the projection $G \to B$ onto the first coordinate
yields a homomorphism $\theta \colon Q \to B$ such that
\[
M = \{(\theta(q), q) : q \in Q\} \leq B \times Q.
\]
Since $M \leq G$, the defining equation of the
fiber product forces $\pi(\theta(q)) = \tau(q)$ for all $q \in Q$.

\emph{Step 3: the coordinate copies of $N$ have a common image under
$\theta$.}

For $1 \leq i \leq p$, let $\varepsilon_i \colon N \to P$ place
$n \in N$ in the $i$th coordinate and $1$ in all others; the image
lies in $P$ because $\pi(n) = 1$. Inside $Q$ we have
$\sigma\varepsilon_i(n)\sigma^{-1} = \varepsilon_{i+1}(n)$, with
indices taken modulo $p$. The order of $\theta(\sigma)$ divides the order $p$ of $\sigma$ as
well as $\card{B}$. Thus, $\theta(\sigma) = 1$ as $\gcd(p, \card{B}) = 1$.
This yields
$\theta\varepsilon_{i+1} = \theta\varepsilon_i$ (composition of maps). Hence
$\psi \coloneqq \theta\varepsilon_i \colon N \to B$ does not depend
on $i$.

\emph{Step 4: $\psi(N') = 1$.}

Tuples supported in distinct coordinates commute, so for all
$m, n \in N$ (using $p \geq 2$),
\[
\psi(m)\psi(n)
= \theta(\varepsilon_1(m))\theta(\varepsilon_2(n))
= \theta\bigl(\varepsilon_1(m)\varepsilon_2(n)\bigr)
= \theta\bigl(\varepsilon_2(n)\varepsilon_1(m)\bigr)
= \psi(n)\psi(m).
\]
Thus $\psi(N)$ is abelian and $N'\le\ker\psi$.

\emph{Step 5: conclusion.}

Let $\delta \colon B \to P$, $b \mapsto (b, \ldots, b)$, be the
diagonal embedding. For $u \in N'$ we have
$\delta(u) = \varepsilon_1(u)\cdots\varepsilon_p(u)$ and therefore
$\theta(\delta(u)) = \psi(u)^{p} = 1$ by Step~4. Hence the
homomorphism $\theta\delta \colon B \to B$ is trivial on $N'$. Moreover,
\begin{equation}\label{eq:piThetaDelta}
\pi(\theta(\delta(b))) = \tau(\delta(b)) = \pi(b)
\end{equation}
for all $b \in B$;
that is, $\theta\delta$ leaves the image in $R$ unchanged.

By the choice of $B$, the subgroup $N/N'$ has a complement in $B/N'$;
let $L \leq B$ be its inverse image, so that $N' \leq L$, $B = NL$
and $N \cap L = N'$. The last two equalities give
$\card{B} = \card{N}\card{L}/\card{N'}$ and hence
$[L : N'] = \card{B}/\card{N} = \card{R}$. Put
$C \coloneqq \theta\delta(L)$. 
Since $N'\le\ker(\theta\delta)$, we obtain $\card{C}\le[L:N']= \card{R}$. 
On the
other hand, 
\[\pi(C) = \pi(L) = \pi(LN)=\pi(B)= R\] 
by \eqref{eq:piThetaDelta}. 
Consequently, $\pi$ restricts to an isomorphism
$C \to R$. In particular $C \cap N = C \cap \ker\pi = 1$, and
$\card{NC} = \card{N}\card{C} = \card{N}\card{R} = \card{B}$
yields $B = NC$. Thus $C$ is a complement of $N$ in $B$,
contradicting the choice of $B$. Therefore $N$ has no complement in
$G$.
\end{proofof}

We remark that \autoref{thm:derived} does not provide a new proof of Gaschütz's theorem because we made use of the original theorem in the proof.

\section{Proof of the main theorem}
\label{sec:finite-extension}

Let $N$ be a group and put $A \coloneqq \aut(N)$.
The inner automorphism map $\gamma \colon N \to A$, $x \mapsto \gamma_x$, is a
homomorphism with image $\inn(N)$ and kernel $\zent(N)$. Since $N'$
is characteristic in $N$, the subgroup $\gamma(N')$ is normal in
$A$. We use the
left-action convention for $\hol(N)=N \rtimes A$, so that
\[
(n, \alpha)(m, \beta) = (n\alpha(m), \alpha\beta)
\qquad (n, m \in N,\ \alpha, \beta \in A).
\]
The map $N\to N\rtimes A$, $x \mapsto (x^{-1}, \gamma_x)$ is a monomorphism,
since for $x, y \in N$
\[
(x^{-1}, \gamma_x)(y^{-1}, \gamma_y)
= \bigl(x^{-1}\gamma_x(y^{-1}), \gamma_x\gamma_y\bigr)
= \bigl(y^{-1}x^{-1}, \gamma_{xy}\bigr)
= \bigl((xy)^{-1}, \gamma_{xy}\bigr).
\]
In particular,
\[
D \coloneqq \bigl\{(d^{-1}, \gamma_d) : d \in N'\bigr\}
             \le \hol(N)
\]
is a subgroup isomorphic to $N'$.
It centralizes $N \times 1$, because for $n \in N$ the
products $(n, 1)(d^{-1}, \gamma_d)$ and $(d^{-1}, \gamma_d)(n, 1)$ are
both equal to $(nd^{-1}, \gamma_d)$. Conjugation by $(1, \alpha)$ with
$\alpha \in A$ sends $(d^{-1}, \gamma_d)$ to
$(\alpha(d)^{-1}, \gamma_{\alpha(d)})$ with $\alpha(d) \in N'$. Since $N \rtimes A$ is generated by
$N \times 1$ and $1 \times A$, we conclude that
$D \unlhd \hol(N)$.

\begin{theorem}\label{thm:universal}
Let $N$ be a group satisfying $\zent(N) \cap N' = 1$. Define
\[
U_N \coloneqq \hol(N)/D,
\qquad
T_N \coloneqq A/\gamma(N').
\]
The following holds:
\begin{enumerate}
\item $N$ embeds as a normal subgroup
of $U_N$.

\item\label{item:piMap} There is an epimorphism $\pi \colon U_N \to T_N$ with kernel
$N$.

\item\label{item:Ncomplement} $N/N'$ has a complement in $U_N/N'$.

\item\label{item:universal-fiber} If $N \unlhd B$ and $N/N'$ has a
complement in $B/N'$, then $B$ is a fiber product of $\pi$ and some
homomorphism $B/N \to T_N$. Here, $N$ is identified with $N\times 1$. 
\end{enumerate}
\end{theorem}

\begin{proofof}\hfill
\begin{enumerate}[(i)]
\item An element $(d^{-1}, \gamma_d)$ of $D$ belongs to
$N \times 1$ if and only if $\gamma_d = 1$, that is, if and only if
$d \in \zent(N) \cap N' = 1$. Hence
$(N \times 1) \cap D = 1$. We can therefore identify $N$ with $(N\times 1)D/D\unlhd U_N$.

\item Since $(n, 1)(d^{-1}, \gamma_d) = (nd^{-1}, \gamma_d)$ for $n\in N$ and $d\in N'$, we have
$(N \times 1)D = N \rtimes \gamma(N')$. This yields a natural homomorphism
\[\pi\colon U_N\to U_N/N\cong (N\rtimes A)/(N\rtimes \gamma(N'))\cong T_N\]
with kernel $N$.

\item As in \ref{item:piMap}, $(N' \times 1)D = N' \rtimes \gamma(N')$. 
Since $\gamma(N')\le\inn(N)$ acts trivially on $N/N'$, $T_N$ acts on $N/N'$. 
This yields a natural isomorphism
\[
U_N/N'\cong (N\rtimes A)/(N'\rtimes\gamma(N'))\cong (N/N')\rtimes T_N.
\]
In particular, $N/N'$ has a complement in $U_N/N'$. 

\item
Now let $N \unlhd B$, and let $L/N'$ be a complement of
$N/N'$ in $B/N'$, i.e. $B=NL$ and $N \cap L = N'$. Since
$N \unlhd B$, conjugation extends $\gamma$ to a homomorphism
$\gamma \colon B \to A$, $\gamma_b(x) \coloneqq bxb^{-1}$ for $x \in N$.
Letting $L$ act on $N$ through $\gamma$, the multiplication map
$\mu \colon N \rtimes L \to B$, $(n, l) \mapsto nl$, is a
homomorphism:
\[
\mu\bigl((n, l)(m, k)\bigr)
= \mu\bigl((nlml^{-1}, lk)\bigr)
= nlmk=\mu(n,l)\mu(m,k).
\]
It is surjective because $B = NL$. For $(n, l) \in \ker\mu$ we have
$l = n^{-1} \in N \cap L = N'$, so
$\ker\mu = \{(d^{-1}, d) : d \in N'\}$. On the other hand,
$N \rtimes L \to N \rtimes A$, $(n, l) \mapsto (n, \gamma_l)$ is a homomorphism
 mapping $\ker\mu$ into $D$,
by the very definition of $D$. Hence the composition
$N \rtimes L \to N \rtimes A \to U_N$ is trivial on $\ker\mu$ and factors
through $\mu$, yielding a homomorphism $f \colon B \to U_N$ with
$f(nl) = (n, \gamma_l)D$ for $n \in N$ and $l \in L$. In
particular, $f(n) = (n, 1)D$ for $n \in N$, so $f$
restricts on $N$ to the identification fixed above. 
The epimorphism $\pi \colon U_N\to U_N/N \cong T_N$ from \ref{item:piMap} induces a homomorphism $\overline{f} \colon B/N \to T_N$ with
$\overline{f}(bN) = \pi(f(b))$.

Consider the fiber product
\[
P \coloneqq U_N \times_{T_N} (B/N)
= \{(u, x) \in U_N \times (B/N) : \pi(u) = \overline{f}(x)\}
\]
and the homomorphism $g \colon B \to P$, $b \mapsto (f(b), bN)$,
which lands in $P$ by the definition of $\overline{f}$. If
$g(b) = 1$, then $b \in N$ and $f(b) = 1$, so $b = 1$ because $f$ is
injective on $N$; hence $g$ is injective. Since $\pi$ is surjective
with kernel $N$, $\card{P} = \card{N}\card{B/N} = \card{B}$.
Hence $g$ is
an isomorphism carrying $N$ onto $N \times 1$, the kernel of the
projection onto $B/N$. \qedhere
\end{enumerate}
\end{proofof}

We can now prove the main theorem stated in the introduction.

\begin{proofofref}{thm:main}
By \cite[Theorem~12]{Sambale2023}, Gaschütz's theorem fails for $N$
if $\zent(N) \cap N' \neq 1$. Assume therefore that
$\zent(N) \cap N' = 1$, so that \autoref{thm:universal} applies.

Suppose first that $N$ has a complement $K$ in $U_N$.
Let $N \unlhd B$ such that $N/N'$ has a complement in $B/N'$. By
\autoref{thm:universal}\ref{item:universal-fiber} we may assume that
$B =P\coloneqq U_N \times_{T_N} (B/N)$ and $N = N \times 1$. The subgroup
\[
M \coloneqq P \cap \bigl(K \times (B/N)\bigr)\le P
\]
intersects $N \times 1$ trivially, because $K \cap N = 1$. Moreover,
let $(u, x) \in P$ and write $u = nk$ with $n \in N$ and $k \in K$.
Since $\ker\pi=N$, we get $\pi(k) = \pi(u)$ and therefore
$(k, x) \in P$, that is, $(k, x) \in M$. Now
$(u, x) = (n, 1)(k, x) \in (N \times 1)M$. Hence
$M$ is a complement of $N$ in $P$. This verifies
\autoref{item:derived-splitting} of \autoref{thm:derived}, so
$\gasch{N}$ holds.

Suppose conversely that $N$ has no complement in $U_N$. By
\autoref{thm:universal}\ref{item:Ncomplement}, $N/N'$ has a
complement in $U_N/N'$, so \autoref{item:derived-splitting} of
\autoref{thm:derived} fails for $B = U_N$. By that theorem,
$\gasch{N}$ fails.
\end{proofofref}

For an abelian group $N$, we have $N' = 1$, so
$U_N = \hol(N)$ splits over $N$. At the other extreme, suppose that
$N$ is perfect and $\zent(N) = 1$. The map
$\hol(N) \to A$, $(n, \alpha) \mapsto \gamma_n\alpha$ is a homomorphism because
$\gamma_{\alpha(n)} = \alpha\gamma_n\alpha^{-1}$. It is clearly surjective with kernel $D$. 
It therefore induces an isomorphism
$U_N \cong A$ carrying $N$ onto $\inn(N)$. Thus
\autoref{thm:main} includes both Gaschütz's theorem and Sambale's
characterization for perfect groups
\cite[Theorem~15]{Sambale2023}.
The following corollary unifies both conditions.

\begin{corollary}\label{cor:automorphism-complement}
Let $N$ be a group satisfying $\zent(N) \cap N' = 1$. If
$\gamma(N')$ has a complement in $A$, then Gaschütz's theorem holds
for $N$.
\end{corollary}

\begin{proofof}
Let $C$ be a complement of $\gamma(N')$ in $A$, and let
\[\overline{C} \coloneqq (1\times C)D/D\le U_N.\] 
An element $(1, c)D$ lies in the image of $N$ if
and only if $(1, c) \in (N \times 1)D = N \rtimes \gamma(N')$,
that is, if and only if $c \in \gamma(N')$. Since $C \cap \gamma(N') = 1$,
we get $\overline{C} \cap N = 1$. Moreover, the epimorphism
$\pi \colon U_N \to T_N$ of
\autoref{thm:universal}\ref{item:piMap} sends 
$\overline{C}$ to $C\gamma(N')/\gamma(N')=T_N$. 
As
$\ker \pi = N$, this gives $N\overline{C} = U_N$. Hence $\overline{C}$
is a complement of $N$ in $U_N$, and Gaschütz's theorem holds for
$N$ by \autoref{thm:main}.
\end{proofof}

The Frobenius group $N\cong (C_3\times C_3)\rtimes Q_8$ shows that the condition in Corollary~\ref{cor:automorphism-complement} is not necessary for $\gasch{N}$.

\section{Solution of a problem of Sambale}\label{sec:problem20}

In \cite{Sambale2023}, the validity of Gaschütz's theorem was
decided for every group of order less than $144$, and the first open
case 
\[N \coloneqq (C_3^2 \rtimes Q_8) \times C_2= \mathtt{SmallGroup}(144,187)\] 
was posed as a problem.
A straightforward computation with GAP~\cite{GAP}\footnote{See the file attached to this pdf.} shows that $N$ has no complement in $U_N \coloneqq \hol(N)/D$ (a group of order $13{,}824$).
A witness $G$ of the failure of $\gasch{N}$ can be constructed via \autoref{thm:derived} with $B=U_N$ and $p=5$. This yields 
\[
\card{G} = p \cdot \card{U_N} \cdot \card{N}^{p}
         = 5 \cdot 13{,}824 \cdot 144^{5}
         = 2^{29} \cdot 3^{13} \cdot 5
         \approx 4.3 \cdot 10^{15}.
\]
The starting point of this paper was a smaller witness of order $2^{25}\cdot 3^2\cdot 5$, found by GPT.

\section*{Acknowledgment}

In preparing this paper the authors made use of the large language
models Claude Opus 5 (Anthropic) and GPT-5.6 Sol (OpenAI).
All definitions, statements and proofs were checked by the authors,
who take full responsibility for the contents of this paper.


\begin{thebibliography}{9}

\bibitem{GAP}
The GAP~Group, \emph{GAP -- Groups, Algorithms, and Programming,
Version 4.16.0}, \url{https://www.gap-system.org}.

\bibitem{Gaschuetz1952}
W.~Gaschütz, \emph{Zur Erweiterungstheorie der endlichen Gruppen},
J. Reine Angew. Math. \textbf{190} (1952), 93--107,
\href{https://doi.org/10.1515/crll.1952.190.93}{doi:10.1515/crll.1952.190.93}.

\bibitem{Sambale2023}
B.~Sambale, \emph{On the converse of Gaschütz' complement theorem},
J. Group Theory \textbf{26} (2023), 931--949,
\href{https://doi.org/10.1515/jgth-2022-0178}{doi:10.1515/jgth-2022-0178}.

\end{thebibliography}
\end{document}